\documentclass[titlepage,11pt]{article}
\usepackage{latexsym}
\usepackage{amsfonts}
\usepackage{amsmath}
\usepackage{mathtools}
\usepackage{tikz}
\usepackage{float}
\usetikzlibrary{decorations.pathreplacing,calligraphy}
\usetikzlibrary{positioning,arrows.meta}
\usepackage{circuitikz}
\usepackage{comment}
\usetikzlibrary{shapes.geometric}
\newcommand\blackslug{\hbox{\hskip 1pt \vrule width 4pt height 8pt depth 1.5pt
        \hskip 1pt}}
\newcommand\bbox{\hfill \quad \blackslug \bigbreak}

\def\DD{\hbox{-}}

\def\LL{,\ldots,}

\title{Distant digraph domination}
\author{
Tung Nguyen\thanks{Supported by AFOSR grant
FA9550-22-1-0234, and by NSF grant DMS-2154169, and a Porter Ogden Jacobus Fellowship.}\\
Princeton University,\\ Princeton, NJ 08544, USA
\and
Alex Scott\thanks{Supported by EPSRC grant EP/X013642/1}\\
University of Oxford, \\
Oxford, UK
\and
Paul Seymour\thanks{Supported by AFOSR grant
FA9550-22-1-0234, and by NSF grant DMS-2154169.}\\
Princeton University,\\ Princeton, NJ 08544, USA}

\date{August 4, 2024; revised \today}

\newtheorem{thm}{}[section]

\newcommand{\Proof}{\noindent{\bf Proof.}\ \ }

\begin{document}
\maketitle
\begin{abstract}
A {\em $k$-kernel} in a digraph $G$ is a stable set $X$ of vertices such that every vertex of $G$ can be joined from $X$ by a directed
path of length at most $k$.  We prove three results about $k$-kernels.

First, it was conjectured by Erd\H{o}s and Sz\'ekely in 1976 that every digraph $G$ with no source has a 2-kernel $|K|$ with 
$|K|\le |G|/2$. We prove this conjecture when $G$ is a ``split digraph'' (that is, its vertex set can be partitioned into a 
tournament and a stable set), improving a result of Langlois et al., who proved that every split digraph $G$ with no source has a 
2-kernel of size at most $2|G|/3$. 

Second, the Erd\H{o}s-Sz\'ekely conjecture implies that in every digraph $G$ there is a 2-kernel $K$ such that the union of $K$ and 
its out-neighbours has size at least $|G|/2$.
We prove that this is true if $V(G)$ can be partitioned into a tournament and an acyclic set.

Third, 
in a recent paper, Spiro asked whether, for all $k\ge 3$, every strongly-connected digraph $G$ has a $k$-kernel of size at 
most about $|G|/(k+1)$. This remains open, but we prove that there is one of size at most about $|G|/(k-1)$. 

\end{abstract}

\section{Introduction}
A {\em digraph} is a finite directed graph with no loops or parallel edges (it may have directed cycles of length two). 
If $G$ is a digraph, $X\subseteq V(G)$ is {\em stable} if there is no edge with both ends in $X$.
In a digraph $G$, if $X,Y\subseteq V(G)$, we say $X$ {\em $k$-covers} $Y$ if for each $y\in Y$, 
there exists $x\in X$ and a directed path of length at most $k$ from $x$ to $y$. (If $X$ is a singleton $\{x\}$ we write $x$ for $\{x\}$ here,
and the same for $Y$.)
A {\em $k$-kernel} in a digraph $G$ is a stable set $X$ of vertices that $k$-covers $V(G)$.
\footnote{In some papers a $k$-kernel is defined with edges reversed: every vertex of $G$ is joined {\em to} $X$ by a short directed path.}

There are many interesting open questions about $k$-kernels; for instance, not every digraph has a 1-kernel, but every digraph has a 2-kernel~\cite{chvatal}, and
the following was conjectured by P. L. Erd\H{o}s and L. A. Sz\'ekely~\cite{erdos}
in 1976
(and remains open):
\begin{thm}\label{smallconj}
{\bf The small quasi-kernel conjecture:} Every digraph $G$ with no source has a 2-kernel of size at most $|G|/2$. 
\end{thm}
(A {\em source} is a vertex with 
in-degree zero.) There is a survey on this conjecture in~\cite{erdossurvey}, and the best bound on this seems to be a result of Spiro~\cite{spiro}, that every digraph $G$ with no source has a $2$-kernel of size at 
most $|G|-\frac14 (|G|\log |G|)^{1/2}$, which is of course very far from the conjecture. 

It is enough to prove \ref{smallconj} for {\em oriented graphs}, that is, digraphs with no directed cycle of length two; 
because deleting an edge from such a cycle makes the problem harder. (Unless this deletion makes a source; but if neither edge will work, 
delete both vertices and all their out-neighbours.) If $G$ is a counterexample to \ref{smallconj}, then, since it has a 2-kernel $S$ say,
it follows that $|S|>|G|/2$; and a natural special case is when $G\setminus S$ is a tournament. 
Let us say $G$ is a {\em split digraph} if $G$ is an oriented graph and its vertex set admits a partition into a stable set and a 
tournament. Ai, Gerke, Gutin, Yeo and Zhou~\cite{ai} proved that \ref{smallconj} holds for split graphs in which all edges between the 
tournament and the stable set are directed towards the stable set.
Langlois, Meunier, Rizzi, Vialette and Zhou~\cite{langlois} proved that every split digraph with no sources admits a 2-kernel of size 
at most $2|G|/3$. In section \ref{sec:splitdi}, we strengthen this:
\begin{thm}\label{mainthm1}
Every split digraph $G$ with no sources admits a 2-kernel $K$ with $|K|\le |G|/2$.
\end{thm}

Our second result concerns a problem of Spiro~\cite{spiro}, who observed that \ref{smallconj} implies:
\begin{thm}\label{largeconj}
{\bf Conjecture:} In every digraph $G$, there is a 2-kernel $K$ such that at least half the vertices of $G$ belong to $K$ or have an
in-neighbour in $K$.
\end{thm}
We discuss this in section \ref{sec:large}, and prove that it holds for split digraphs, and indeed for digraphs with a vertex set 
that can be partitioned into a tournament and an acyclic subgraph.

Our third result concerns a different problem of Spiro~\cite{spiro}, who asked whether:
\begin{thm}\label{openconj}
{\bf Conjecture:} For all integers $k\ge 3$, every strongly-connected digraph $G$ has a $k$-kernel of size at most $|G|/(k+1)+O_k(1)$.
\end{thm}
It seems that the best known bound in this case is due to Spiro, in the same paper, who proved that under the hypotheses of 
\ref{openconj}, there is a $k$-kernel of size at most about $|G|/\log k$. Our third result is that there is one of size at most 
$|G|/(k-1)+O_k(1)$. This as a consequence of \ref{arbthm} below.

Let $T$ be a subdigraph with underlying graph a tree, such that for some vertex $r$ of $T$, every edge of $T$ is directed away 
from $r$ in the natural sense. We call $T$ an {\em arborescence}, and $r$ is its {\em root}. Every strongly-connected digraph has 
a subdigraph that is a 
spanning arborescence ({\em spanning} means that the arborescence contains all vertices of the digraph). 
In section \ref{sec:kkernels} we will prove:
\begin{thm}\label{arbthm}
For all integers $k\ge 2$, every digraph $G$ with $|G|>1$ and with a spanning arborescence has a $k$-kernel of size at most $1+(|G|-2)/(k-1)$.
\end{thm}

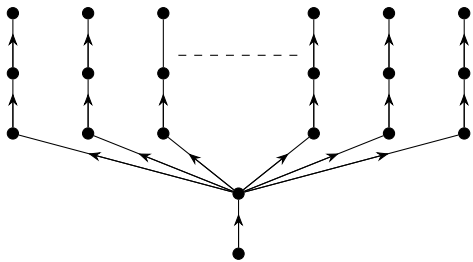
\begin{figure}[ht]
\centering

\begin{tikzpicture}[auto=left]
\tikzstyle{every node}=[inner sep=1.5pt, fill=black,circle,draw]

\def\s{.8}
\node (x) at (0,0) {};
\node (y) at (0,-\s) {};
\node (a1) at (-3,\s) {};
\node (a2) at (-2,\s) {};
\node (a3) at (-1,\s) {};
\node (a4) at (1,\s) {};
\node (a5) at (2,\s) {};
\node (a6) at (3,\s) {};
\node (b1) at (-3,2*\s) {};
\node (b2) at (-2,2*\s) {};
\node (b3) at (-1,2*\s) {};
\node (b4) at (1,2*\s) {};
\node (b5) at (2,2*\s) {};
\node (b6) at (3,2*\s) {};
\node (c1) at (-3,3*\s) {};
\node (c2) at (-2,3*\s) {};
\node (c3) at (-1,3*\s) {};
\node (c4) at (1,3*\s) {};
\node (c5) at (2,3*\s) {};
\node (c6) at (3,3*\s) {};

\tikzstyle{every node}=[];
\foreach \from/\to in {y/x,x/a1,x/a2,x/a3,x/a4,x/a5,x/a6,a1/b1,a2/b2,a3/b3,a4/b4,a5/b5,a6/b6, b1/c1,b2/c2,b3/c3,b4/c4,b5/c5,b6/c6}
\draw [] (\from)--(\to);
\draw [-Stealth] (y) -- (0,-\s/3);
\draw [-Stealth] (x) -- (-2,2*\s/3);
\draw [-Stealth] (x) -- (-4/3,2*\s/3);
\draw [-Stealth] (x) -- (-2/3,2*\s/3);
\draw [-Stealth] (x) -- (2/3,2*\s/3);
\draw [-Stealth] (x) -- (4/3,2*\s/3);
\draw [-Stealth] (x) -- (2,2*\s/3);
\draw [-Stealth] (a1) -- (-3,5*\s/3);
\draw [-Stealth] (a2) -- (-2,5*\s/3);
\draw [-Stealth] (a3) -- (-1,5*\s/3);
\draw [-Stealth] (a4) -- (1,5*\s/3);
\draw [-Stealth] (a5) -- (2, 5*\s/3);
\draw [-Stealth] (a6) -- (3,5*\s/3);
\draw [-Stealth] (b1) -- (-3,8*\s/3);
\draw [-Stealth] (b2) -- (-2,8*\s/3);
draw [-Stealth] (b3) -- (-1,8*\s/3);
\draw [-Stealth] (b4) -- (1,8*\s/3);
\draw [-Stealth] (b5) -- (2,8*\s/3);
\draw [-Stealth] (b6) -- (3,8*\s/3);
\draw[dashed] (-.8,2.3*\s) -- (.8,2.3*\s);

\end{tikzpicture}

\caption{All 3-kernels have size $\ge 1+(|G|-2)/3$. For $k>3$ make the vertical paths longer. } \label{fig:3kernels}
\end{figure}

This follows easily from a result about acyclic digraphs ({\em acyclic} means there is no directed cycle):
\begin{thm}\label{acyclic}
For every integer $k\ge 1$, if $G$ is an acyclic digraph with $|G|\ge 2$ and with only one source, then $G$ has a $k$-kernel of size at most $1+(|G|-2)/k$.
\end{thm}
This result is tight, as can be seen from the digraph shown in figure \ref{fig:3kernels}.


\section{Split digraphs}\label{sec:splitdi}

If $G$ is a digraph, we use $G[X]$ to denote the subdigraph induced on $X\subseteq V(G)$. We say ``$u$ is adjacent to $v$'' to mean that
$u$ is an in-neighbour of $v$, and ``adjacent from'' to mean it is an out-neighbour.  A {\em neighbour} of $v$ means a vertex that is 
either an in-neighbour or an out-neighbour of $v$. We sometimes use ``$G$-in-neighbour'' to mean 
``in-neighbour in the digraph $G$'', and so on (this is helpful because we 
sometimes work with different digraphs that have the same vertex set.)
For a vertex $v$ of a digraph $G$, $N^+_G(v)$ denotes the set of all out-neighbours of $v$, and $N^-_G(v)$ is its set of in-neighbours.
A {\em split} in an oriented graph $G$ is a pair $(S,T)$, where $S\cup T=V(G)$, $S\cap T=\emptyset$, $S$ is a stable set, and $G[T]$ is 
a tournament. (We will often write $T$ for $G[T]$.) 

In this section we prove \ref{mainthm1}, but it is convenient to prove a slightly stronger statement, that the same conclusion holds 
just assuming that no vertex in $S$ is a source. Now there is a difficulty, because this is false for the 1-vertex digraph with 
$S=\emptyset$,
but this is the only exception. We will prove:
\begin{thm}\label{mainthm2}
Let $(S,T)$ be a split of an oriented graph $G$, such that $S\ne \emptyset$ and no vertex in $S$ is a source. Then there is a 2-kernel $K$
with $|K|\le |G|/2$.
\end{thm}

For the proof, we begin with some lemmas.
A 2-kernel $K$ is {\em strong} if for every vertex $v\in T$, either there is a vertex in 
$K$ that 1-covers $v$, or a vertex in $K\cap T$ that 2-covers $v$.
(We do not know whether \ref{mainthm1} remains true if we ask for a strong 2-kernel of size at most $|G|/2$.)
If $v\in T$, we say $s\in S$ is a {\em problem} for $v$ if $v$ is adjacent from $s$, and $v$ does not 2-cover $s$, 
and no non-neighbour of $v$ in $S$ 2-covers $s$. If $v$ has a problem, then $v$ is contained in no 2-kernel.

\begin{thm}\label{goodvert}
Let $G,T,S$ be as above, and let $v\in V(T)$. If $v$ is contained in no strong 2-kernel, then there exists $w\in V(T)\setminus \{v\}$,
adjacent to $v$, such that $N^-_G(w)\subseteq N^-_G(v)$; and either $w\in S$ and $w$ is a problem for $v$, or $w\in T$.
\end{thm}
\Proof
Since the set consisting of $v$ and all non-neighbours of $v$ in $S$ is not a strong 2-kernel, there exists $w\in V(G)\setminus \{v\}$ such that $v$ does not 2-cover $w$,
and either $w\in T$ and no non-neighbour of $v$ in $S$ 1-covers $w$, or $w\in S$ and no non-neighbour of $v$ in $S$ 2-covers $w$. 
In the first case, since $v$ does 
not 2-cover $w$, $N^-_G(w)\cap T\subseteq N^-_G(v)$.
If $s\in N^-_G(w)\cap S$, then since no non-neighbour of $v$ in $S$ 1-covers $w$, it follows that
$s\in N^+_G(v)\cup N^-_G(v)$; and since $v$ does not 2-cover $w$, $s\notin N^+_G(v)$, and so $s\in  N^-_G(v)$.
This proves that $N^-_G(w)\subseteq N^-_G(v)$ as required.
In the second case, $w$ is a problem for $v$. Moreover, every in-neighbour of $w$ is an in-neighbour of $v$: because if $u\in T$
is adjacent to $w$, then $u$ is not adjacent from $v$ since $v$ does not 2-cover $w$, and so $u$ is adjacent to $v$.
Hence, again, $N^-_G(w)\subseteq N^-_G(v)$. This proves \ref{goodvert}.~\bbox

\begin{thm}\label{nondom}
Let $G,T,S$ be as above, and suppose that $G, S,T$ form a smallest counterexample to \ref{mainthm2}.
Suppose also that $v\in V(T)$ is contained in no strong 2-kernel, and let $w$ be as in \ref{goodvert}. If $w\in T$, then there is no problem 
for $w$.
\end{thm}
\Proof Suppose that $w\in T$, and $s\in S$ is a problem for $w$. Let $A=N_G^+(v)$.
Since $N^-_G(w)\subseteq N^-_G(v)$, no 
vertex in $A$ is adjacent to $w$, and in particular $s\notin A$.
Make a digraph $G'$ from $G$ by deleting $v$ and making $w$ complete to $A$. So $G'$ has no sources.
\\
\\
(1) {\em $N^-_{G'}(w) \subseteq N^-_G(v)$.}
\\
\\
Let $u\in N^-_{G'}(w)$. So $u\notin A$, and so $u\in N^-_G(w)\subseteq N^-_G(v)$. This proves (1).

\bigskip

Let $K$ be a 2-kernel of $G'$. We will show that $K$ is also a 2-kernel of $G$.
Certainly it is stable in $G$.
\\
\\
(2) {\em $w\notin K$.}
\\
\\
Suppose that $w\in K$.
Then $s\notin K$, so there is a directed path $P$ of $G'$, of length one or two, from some $x\in K$ to $s$. Since $s$ is a problem 
for $w$ in $G$, some edge of $P$ is not an edge of $G$, which is impossible since $s\notin A$.
This proves (2).

\bigskip

So $w\notin K$. Since $K$ 2-covers $w$ in $G'$, (1) implies that $K$ 2-covers $v$ in $G$, and 1-covers $v$ in $G$ if it 1-covers $w$ in $G'$. 
Let $a\in A$. We must show that $K$ 2-covers $a$ in $G$.
If $a\in K$ this is true, so we assume
there is a directed path $P$ of $G'$ of length one or two, from some $x\in K$ to $a$. If $P$ is a path of $G$ then $K$ 2-covers $a$ 
in $G$, so we may assume that
the last edge of $P$ is an edge of $G'$ not in $G$. But $w\notin K$ and $x\in K$, so $w\ne x$, and therefore $P$ has length two
with middle vertex $w$. By (1),
$x\DD v\DD a$ is a path of $G$, so $K$ 2-covers $a$ in $G$.

This proves that every 2-kernel of $G'$ is a 2-kernel of $G$. Since $G,S,T$ form a smallest counterexample to \ref{mainthm2},
and $G'$ has fewer vertices than $G$, and $(S,T\setminus \{v\})$ is a split for $G'$, with $S\ne \emptyset$, and no vertex in $S$ is a source in $G'$, it follows that $G'$ has a 2-kernel of size at most $|G'|/2 $; but this is also a 2-kernel for $G$,
which is impossible. This proves that there is no problem for $v$, and so proves \ref{nondom}.~\bbox

Now we prove the main theorem, which we restate:
\begin{thm}\label{mainthm3}
Let $(S,T)$ be a split of an oriented graph $G$, such that $S\ne \emptyset$ and no vertex in $S$ is a source. Then there is a 2-kernel $K$
with $|K|\le |G|/2$.
\end{thm}
\Proof
We may assume that $G,S,T$ form a smallest counterexample. 
Let $B$ be the set of all vertices in $T$ with problems.
For each $b\in B$, select a problem $z_b$ for $b$, and let $Z$ be the set $\{z_b:b\in B\}$.
Let $Q$ be the set of all $q\in S\setminus Z$ with $N^-_G(q)\subseteq B$.
For each $q\in Q$, it has an in-neighbour in $B$, since it is not a source;
select one such in-neighbour $b_q$.
Similarly, for each $s\in S\setminus (Q\cup Z)$, choose some $t_s\in T\setminus B$ adjacent to $s$.

For each $z\in Z$, let $\Phi(z)$ be the set of $q\in Q$ such that $z=z_{b_q}$.
For each $t\in T\setminus B$, let $\Phi(t)$ be the union of $\{t\}$ and the set of $s\in S\setminus (Q\cup Z)$ such that $t=t_s$.
Thus, the sets $\Phi(v)\;(v\in V(H))$ are pairwise disjoint and have union $V(G)\setminus (B\cup Z)$. Some of the sets $\Phi(z)\;(z\in Z)$
may be empty.
\begin{figure}[ht]
\centering

\begin{tikzpicture}[auto=left]
\tikzstyle{every node}=[inner sep=1.5pt, fill=black,circle,draw]
\draw[rounded corners] (-4,-.7) rectangle (5,0);
\draw[rounded corners] (-7,1) rectangle (7,2.25);
\draw[-] (2,-.7) -- (2,0);
\draw[-] (0,1) -- (0,2.25);
\draw[-] (2,1) -- (2,2.25);
\node (z1) at (.5,1.5) {};
\node (z2) at (1.5,1.5) {};
\node (b1) at (2.5,-.35) {};
\node (b2) at (3,-.35) {};
\node (b3) at (3.5,-.35) {};
\node (q1) at (2.5,1.5) {};
\node (q2) at (3.25,1.5) {};
\node (q3) at (4,1.5) {};
\node (q4) at (4.75,1.5) {};
\node (q5) at (5.5,1.5) {};
\node (q6) at (6.25,1.5) {};
\node (t1) at (-3,-.35) {};
\node (t2) at (-2,-.35) {};
\node (t3) at (-1,-.35) {};
\node (s1) at (-6,1.5) {};
\node (s2) at (-5,1.5) {};
\node (s3) at (-4,1.5) {};
\node (s4) at (-3,1.5) {};
\node (s5) at (-2,1.5) {};
\node (s6) at (-1,1.5) {};

\tikzstyle{every node}=[];
\foreach \from/\to in {z1/b1,z2/b2,z2/b3,b1/q1,b1/q2,b2/q3,b2/q4,b3/q5,b3/q6,t1/s1,t1/s2,t2/s3,t2/s4,t3/s5,t3/s6}
\draw [] (\from)--(\to);
\draw [-{Stealth[length=3mm, width=2mm]}] (z1) -- (3/2,1.15/2);
\draw [-{Stealth[length=3mm, width=2mm]}] (z2) -- (4.5/2,1.15/2);
\draw [-{Stealth[length=3mm, width=2mm]}] (z2) -- (5/2,1.15/2);
\draw [-{Stealth[length=3mm, width=2mm]}] (b1) -- (2.5,1.15/2);
\draw [-{Stealth[length=3mm, width=2mm]}] (b1) -- (5.75/2,1.15/2);
\draw [-{Stealth[length=3mm, width=2mm]}] (b2) -- (7/2,1.15/2);
\draw [-{Stealth[length=3mm, width=2mm]}] (b2) -- (7.75/2,1.15/2);
\draw [-{Stealth[length=3mm, width=2mm]}] (b3) -- (9/2,1.15/2);
\draw [-{Stealth[length=3mm, width=2mm]}] (b3) -- (9.75/2,1.15/2);
\draw [-{Stealth[length=3mm, width=2mm]}] (t1) -- (-9/2,1.15/2);
\draw [-{Stealth[length=3mm, width=2mm]}] (t1) -- (-8/2,1.15/2);
\draw [-{Stealth[length=3mm, width=2mm]}] (t2) -- (-6/2,1.15/2);
\draw [-{Stealth[length=3mm, width=2mm]}] (t2) -- (-5/2,1.15/2);
\draw [-{Stealth[length=3mm, width=2mm]}] (t3) -- (-3/2,1.15/2);
\draw [-{Stealth[length=3mm, width=2mm]}] (t3) -- (-2/2,1.15/2);

\node at (4.5,-.35) {$B$};
\node at (.5,-.35) {$T\setminus B$};
\node at (3,1.85) {$Q$};
\node at (.7,1.85) {$Z$};
\node at (-4,1.85) {$S\setminus (Q\cup Z)$};
\node at (5.125,1.85) {$\Phi(z)$};
\node at (1.5,1.7) {$z$};
\node at (-1.6,1.85) {$\Phi(t)$};
\node at (-.8,-.35) {$t$};

\draw (5.125,1.6) ellipse (1.6 and .5);
\draw[rotate around={110:(-1.3,.8)}] (-1.3,.8) ellipse (1.5 and .7);

\end{tikzpicture}

\caption{Definitions of $\Phi(z)$ and $\Phi(t)$. } \label{fig:Phi}
\end{figure}
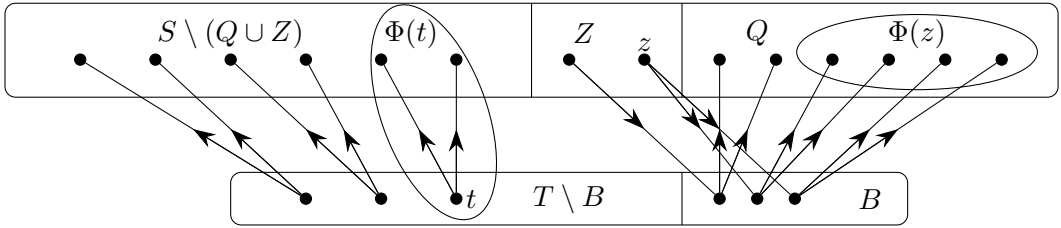

Let $H$ be the oriented graph obtained from $G[(T\setminus B)\cup Z]$
by adding all possible edges from $T\setminus B$ to $Z$; that is, if $t\in T\setminus B$ and $z\in Z$ are nonadjacent in $G$
then we add an edge $tz$. 

For each $v\in V(H)$, let 
$N^0_H(v)$ be the set of vertices that are neither out- nor in-neighbours of $v$ (including $v$ itself). Thus $N^0_H(v)=Z$ if $v\in Z$, 
and $N^0_H(v)=\{v\}$ if $v\in V(T)\setminus B$. Define
$\phi^+(v) = \sum_{u\in N^+_H(v)}|\Phi(u)|$ and define $\phi^-(v), \phi^0(v)$ similarly. We call $\phi^-(v) + \phi^0(v)/2$ the 
{\em score} of $v$. If $V(H)=\emptyset$, then $T\setminus B=\emptyset$ and $B=\emptyset$ (since $Z=\emptyset$); so $T=\emptyset$, 
which implies that $S=\emptyset$ (since there are no sources), a contradiction. So $V(H)\ne \emptyset$. 
We have 
$$\sum_{u\in V(H)} |\Phi(u)| \phi^+(u) = \sum_{uw\in E(H)} |\Phi(u)||\Phi(w)| = \sum_{w\in V(H)} |\Phi(w)|\phi^-(w),$$
and therefore
$$\sum_{u\in V(H)}|\Phi(u)| (\phi^-(u)-\phi^+(u)) = 0.$$

We claim that there exists 
$v\in V(H)$ such that $\phi^+(v)\ge \phi^-(v)$. If $|\Phi(u)| (\phi^-(u)-\phi^+(u)) \ne 0$ for some $u\in V(H)$, then 
$|\Phi(u)| (\phi^-(u)-\phi^+(u)) >0$ for some $u\in V(H)$ and the claim is true. If not, then either $|\Phi(u)|=0$ for each $u\in V(H)$, or 
$\phi^-(u)-\phi^+(u)=0$
for some $u\in V(H)$, and in either case the claim is true.
This proves that there exists $v\in V(H)$ such that $\phi^+(v)\ge \phi^-(v)$. 

Since 
$$\phi^+(v)+\phi^-(v)+\phi^0(v)=|G|-|Z|-|B|\le |G|-2|Z|,$$ 
it follows that
$\phi^-(v) + \phi^0(v)/2\le |G|/2-|Z|$. 
Choose $v\in V(H)$ with score as small as possible (and consequently with score at most $|G|/2-|Z|$).

A vertex in $T$ is {\em pure-up} if it has no in-neighbour in $S$.
The case when $v$ has score exactly $|G|/2-|Z|$ is troublesome, so let us first handle that. 
\\
\\
(1) {\em We may assume that either $v$ has score strictly less than $|G|/2-|Z|$, or $v\in Z$ and $\Phi(v)\ne \emptyset$, or $|\Phi(v)|\ge 2$.}
\\
\\
We assume that $v$ has score exactly $|G|/2-|Z|$. It follows that $|B|=|Z|$, and 
every vertex $u\in V(H)$ has score at least $|G|/2-|Z|$, and so 
satisfies $\phi^+(u)\le \phi^-(u)$. 
But 
$$\sum_{u\in V(H)}|\Phi(u)| (\phi^-(u)-\phi^+(u)) = 0.$$
It follows that for every $u\in V(H)$, $|\Phi(u)|(\phi^-(u)-\phi^+(u))=0$, so either $\Phi(u)=\emptyset$ (and hence $u\in Z$) or $\phi^+(u)=\phi^-(u)$ (and hence $u$ has 
the same score as $v$). In particular, if $\Phi(u)\ne\emptyset$ for some $u\in Z$, then we may replace $v$ by $u$ and the claim holds.
Similarly, if some $u\in T\setminus B$
satisfies $|\Phi(u)|\ge 2$, we can replace $v$ by $u$. So we may assume that $\Phi(u)=\emptyset$ for all $u\in Z$ 
(and hence $Q=\emptyset$), and $\Phi(u)=\{u\}$ for each $u\in T\setminus B$ (and hence $S\setminus (Q\cup Z)=\emptyset$). Consequently,
$S=Z$.
Since $|Z|\le |G|/2$ (because $|Z|=|B|$), we may assume that there exists $p_0\in T$ not 2-covered by $Z$. Thus $p_0$ is pure-up,
and so $P\ne \emptyset$, 
where $P$ is the
set of pure-up vertices. 
Choose $p\in P$ that 2-covers $P$. (Any vertex of maximum out-degree in $T[P]$ has this property.) Let $Z'$ be the set of vertices in 
$Z$ that are not adjacent from $p$; so $Z'\cup \{p\}$ is stable. We claim it is a 2-kernel. Certainly $Z'\cup \{p\}$ 2-covers 
$Z$; each vertex in $T$ 1-covered by $Z'$ is 2-covered by $p$; every other vertex of $T$ 1-covered by $Z$
is 2-covered by $Z\setminus Z'$; and each vertex of $T$ not 1-covered by $Z$ is in $P$, and hence is 2-covered by $p$. 
So $Z'\cup \{p\}$ is a 2-kernel, and therefore we may assume its size is more than $|G|/2$. Since $|Z|=|B|$, it follows that 
$|T\setminus B|=1$
and hence $T\setminus B=P=\{p\}$, since $P\cap B=\emptyset$; and so $p_0=p$. Since $Z$ 1-covers $B$ and does not 
2-cover $p_0=p$,
it follows that $p$ is adjacent to every vertex in $B$. But then $\{p\}$ is a 2-kernel (because every vertex in $S=Z$
has an in-neighbour, since it is not a source). This proves (1).
\\
\\
(2) {\em If $v\in Z$ then the theorem holds.}
\\
\\
Let $J$ be the set of vertices in $S\setminus Z$ that are 2-covered 
by $v$. (Possibly $J\cap Q\ne \emptyset$.) 
Let $A=S\setminus (J\cup Q\cup Z)$, and $F=(T\setminus B)\setminus N^+_G(v)$. Since  $N^-_H(v)= F$, and therefore the union of the sets 
$\Phi(u)\;(u\in N^-_H(v))$ includes $F\cup A$, it follows that
$\phi^-(v)\ge  |F| + |A|$.
Moreover, 
$$\phi^0(v) = \frac12 \sum_{z\in Z}|\Phi(z)|=|Q|/2.$$
Consequently, the score of $v$ is at least 
$|F| + |A|+ |Q|/2$,
and so the latter is at most $|G|/2-|Z|.$

Choose $X\subseteq S$ minimal such that $A\cup Z\cup X$ 1-covers every vertex of $T$ that is not pure-up.
Thus $|X|\le |F|$, since $Z$ 1-covers $B\cup (T\cap N^+_G(v))$.
Let $K= A\cup Z\cup X$. We claim that $K$ is a 
2-kernel. It certainly 2-covers $S$, since $Z$ 2-covers $Q$, and $A\cup \{v\}$ 2-covers $S\setminus (Q\cup Z)$. 
It 1-covers all vertices in $T$ that are not pure-up, from the choice of $X$.
Suppose it does not 2-cover some $p\in T\setminus B$. Then $p$ is pure-up, so $p\notin B$; and 
$p$ is complete to all vertices in $T$ that are not pure-up, since $X$ 1-covers all such vertices and does not 2-cover $p$. 
Moreover, each vertex in $Z$ is adjacent from $p$ in $H$.
Thus, every $H$-in-neighbour of $p$ is also pure-up, and so is adjacent to $v$ in $H$. Consequently
$$|\Phi(p)|+\sum_{u\in N^-_H(p)}|\Phi(u)| \le \sum_{u\in N^-_H(v)}|\Phi(u)|;$$
and so  $p$ has smaller score than $v$, a contradiction.

So $K$ is a 2-kernel. 
But 
$$|K|\le |X|+|A|+|Z|\le |F|+|A|+|Z|\le |G|/2-|Q|/2.$$
It follows that $|K|\le |G|/2$. This proves (2).

\bigskip

Henceforth we assume that $v\in V(T)\setminus B$ and, by (1),  either $v$ has score strictly less than $|G|/2-|Z|$, or $|\Phi(v)|\ge 2$. 
\\
\\
(3) {\em $v$ extends to a strong 2-kernel.}
\\
\\
Suppose not. By \ref{goodvert},  there exists $t\in T$, adjacent to $v$, such that every $G$-in-neighbour of $t$ is 
a $G$-in-neighbour of $v$, and $t\in T\setminus B$ by \ref{nondom}. A vertex of $H$ is a $G$-in-neighbour of $v$ if and only if it is an $H$-in-neighbour of $v$, and the 
same is true for in-neighbours of $t$; so  every
$H$-in-neighbour of $t$ is
an $H$-in-neighbour of $v$. 
Hence $\phi^-(v)\ge \phi^-(t) + |\Phi(t)|$. Since $\phi^0(v)=|\Phi(v)|$ and $\phi^0(t)=|\Phi(t)|$,
it follows that 
$$\phi^-(v)+\phi^0(v)/2\ge \phi^-(t) + |\Phi(t)| + |\Phi(v)|/2> \phi^-(t)+\phi^0(t)/2,$$
and so the score of $t$ is strictly less than that of $v$,
contradicting the choice of $v$.
This proves (3).

\bigskip

\begin{figure}[H]
\centering

\begin{tikzpicture}[auto=left]
\draw[rounded corners] (0,0)  rectangle (12,3);
\draw[rounded corners] (2,-2)  rectangle (10,-1);

\draw[-] (2,0) -- (2,3);
\draw[-] (4,0) -- (4,3);
\draw[-] (6,0) -- (6,3);
\draw[-] (8,0) -- (8,3);
\draw[-] (10,0) -- (10,3);
\draw[-] (0,1) -- (12,1);
\draw[-] (0,2) -- (12,2);

\draw[-] (4,-2) -- (4,-1);
\draw[-] (6,-2) -- (6,-1);
\draw[-] (8,-2) -- (8,-1);

\tikzstyle{every node}=[inner sep=1.5pt, fill=black,circle,draw]
\node (v) at (-1,-.5) {};
\tikzstyle{every node}=[]

\node at (1,2.5) {$\emptyset$};
\node at (7,2.5) {$\emptyset$};
\node at (9,2.5) {$\emptyset$};
\node at (11,2.5) {$\emptyset$};
\node at (13,2.5) {$S\cap N_G^+(v)$};
\node at (13,1.5) {$S\setminus N_G(v)$};
\node at (13,0.7) {$S\cap N_G^-(v)$};
\node at (12.8,0.3) {$=D$};
\node at (9,3.5) {$Q''$};
\node at (11,3.5) {$Q'$};
\node at (4,3.5) {$J$};
\node at (6,3.5) {$Z$};
\draw [decorate,
    decoration = {calligraphic brace, amplitude = 9}] (2,3) --  (6,3);
\draw [decorate,
    decoration = {calligraphic brace,amplitude = 6,raise = 3}] (4,3) --  (8,3);

\draw[-{Stealth}, line width = 1] (v) to [out=90, in=180, looseness=1] (0,2.5);
\draw[-{Stealth}, line width = 1] (0,.5) to [out=180, in=30, looseness=1] (v);
\draw[-{Stealth}, line width = 1] (v) to [out=330, in=180, looseness=.7] (3,-1.5);
\draw[-{Stealth}, line width = 1] (4.8,-1.5) to [out=200, in=270, looseness=.7] (v);
\draw[-{Stealth}, line width = 1] (2.4,-1.2) to (3,1.5);
\draw[-{Stealth}, line width = 1] (2.8,-1.2) to (3,0.5);
\draw[-{Stealth}, line width = 1] (3.2,-1.2) to (4.7,1.5);
\draw[-{Stealth}, line width = 1] (3.6,-1.2) to (4.7,0.5);

\draw[-{Stealth}, line width = 1] (5.2,.5) to (6.5,-1.5);
\draw[-{Stealth}, line width = 1] (6.2,.5) to (7.5,-1.5);

\draw[-{Stealth}, line width = 1] (5.6,1.5) to (8.4,-1.5);
\draw[-{Stealth}, line width = 1] (5.6,2.14) to (9.0,-1.5);
\draw[-{Stealth}, line width = 1] (6.8,1.5) to (9.6,-1.5);

\draw[-{Stealth}, line width = 1] (6.5,-1.5) to (8.5,1.5);
\draw[-{Stealth}, line width = 1] (6.5,-1.5) to (8.5,0.5);
\draw[-{Stealth}, line width = 1] (7.5,-1.5) to (9.5,1.5);
\draw[-{Stealth}, line width = 1] (7.5,-1.5) to (9.5,.5);

\draw[-{Stealth}, line width = 1] (8.4,-1.5) to (10.4,1.5);
\draw[-{Stealth}, line width = 1] (8.4,-1.5) to (10.4,0.5);
\draw[-{Stealth}, line width = 1] (9,-1.5) to (11,1.5);
\draw[-{Stealth}, line width = 1] (9,-1.5) to (11,0.5);
\draw[-{Stealth}, line width = 1] (9.6,-1.5) to (11.6, 1.5);
\draw[-{Stealth}, line width = 1] (9.6,-1.5) to (11.6, 0.5);

\node[left] at (v) {$v$};
\draw [decorate,
    decoration = {calligraphic brace,mirror, amplitude = 6,raise = 3}] (6.1,-2) --  (10,-2);
\node at (8,-2.5) {$B$};
\draw [decorate,
    decoration = {calligraphic brace,mirror, amplitude = 6,raise = 3}] (2,-2) --  (5.9,-2);
\node at (4,-2.5) {$T\setminus (B\cup \{v\})$};
\node at (5,-1.5) {$F$};

\end{tikzpicture}

\caption{$v$ is adjacent to everything in the top row of boxes, and from everything in the third. Its adjacency to $B$ is not specified in the figure. It has no out-neighbours in $Q$ since $v\notin B$, and so all its out-neighbours in $S$ belong to $J$.} \label{fig:Spart}
\end{figure}
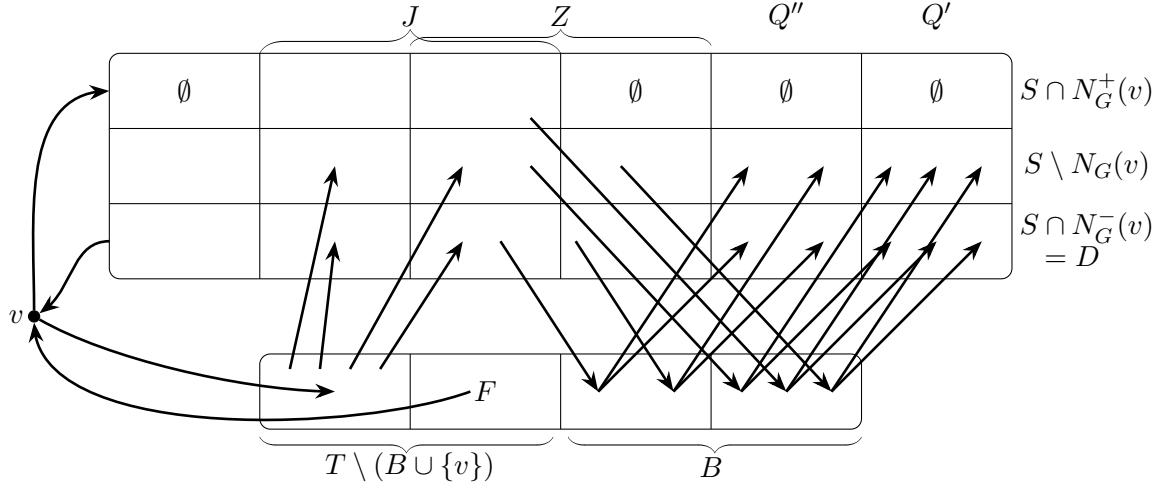

Let $Q'=\bigcup_{z\in Z\setminus N^-(v)} \Phi(z)$, and $Q''=\bigcup_{z\in Z\cap N^-(v)} \Phi(z)$; so $Q''=Q\setminus Q'$.
Let $J$ be the set of vertices in $S\setminus Q$ that are 2-covered by $v$ in $G\setminus B$. 
So, $J,Z$ are both subsets of $S\setminus Q$, but they might intersect each other. $S$ is also partitioned into three subsets,
$S\cap N^+_G(v)$, $S\cap N^-_G(v)$ and $S\setminus N_G(v)$, where we define $N_G(v) =N^+_G(v)\cup N^-_G(v)$.
(See figure \ref{fig:Spart}.) We intend to find a 2-kernel containing $v$ of size at most $|G|/2$,
but we must be careful only to add vertices in $S\setminus N_G(v)$, to keep the set stable.

Let $D=N^-_G(v)\cap S$, and $F=(T\setminus B)\cap N^-_G(v)$.
Thus 
$$N^-_H(v)= F\cup (Z\cap D).$$
The union of the sets $\Phi(t)\;(t\in F)$ includes 
$F \cup (S\setminus (Q\cup J\cup Z))$,
and $\bigcup_{z\in Z\cap D}\Phi(z)= Q''$.
Consequently 
$$\phi^-(v)\ge |F|+ |S\setminus (Q\cup J\cup Z)|+ |Q''|,$$
and so the score of $v$ is at least      
$$|F|+|S\setminus (Q\cup J\cup Z)|+ |Q''|+ \phi^0(v)/2.$$
Since $\phi^0(v)\ge 1$, and either $\phi^0(v)\ge 2$ or the score of $v$ is strictly less than $|G|/2-|Z|$, 
it follows that 
$$|F|+|S\setminus (Q\cup J\cup Z)|+ |Q''|+ 1+ |Z|\le |G|/2.$$

Since $v$ extends to a strong 2-kernel, for each $u\in T\setminus B$ that is not 2-covered by $v$, there is an 
in-neighbour of $u$ in $S\setminus N_G(v)$;
choose $X\subseteq S\setminus N_G(v)$ minimal 1-covering each vertex in $F$ that is not 2-covered by $v$. 
Thus $|X|\le |F|$. For each 
$u\in D$, since $v$ extends to a 2-kernel, there exists $t\in S\setminus N_G(v)$ that 2-covers $u$;
let $Y\subseteq S\setminus N_G(v)$ be minimal 2-covering $D\setminus (J\cup Q')$. Thus $|Y|\le |D\setminus (J\cup Q')|$.

Let 
$$K=\{v\}\cup (Z\setminus D) \cup (S\setminus (Q\cup J\cup Z\cup D))\cup 
X\cup Y\cup (Q''\setminus D).$$
We claim that $K$ is a 2-kernel. Certainly it is stable. 
\\
\\
(4) {\em $K$ 2-covers $S$.}
\\
\\
Let $s\in S$, and assume first that $s\notin Q$. If $s\in J$ then $v$ 2-covers $s$; 
if $s\in D\setminus J$ then $Y$ 2-covers $s$; if
if $s\in Z\setminus (J\cup D)$ then $s\in K$; and if $s\notin Z\cup J\cup D$ then $s\in K$. So in this case $K$
2-covers $s$. Next assume that $s\in Q$. So $s\notin J\cup N^+_G(v)$. If $s\in Q''\setminus D$ then $s\in K$, and  if 
$s\in Q''\cap D$ then $Y$ 2-covers $s$, so we assume that $s\in Q'$, and so $z_{b_s}\in Z\setminus D$. If 
$z_{b_s}\notin N^+_G(v)$
then $z_{b_s}\in K$ and so $K$ 2-covers $s$, so we assume that $z_{b_s}\in N^+_G(v)$.
Then $b_s$ is adjacent from $v$ (because $b_s$
does not 2-cover $z_{b_s}$ since $z_{b_s}$ is a problem for $b_s$)
and so $K$ 2-covers $s$. This proves (4).
\\
\\
(5) {\em $K$ 2-covers $T$, and hence $K$ is a 2-kernel.}
\\
\\
Let $t\in T$. We may assume that $t\in N^-_G(v)$. If 
$t\in T\setminus B$ then $t\in F$ and $X$ 1-covers $t$, so we assume that $t\in B$. If $z_t\notin N_G(v)$ then $z_t\in K$
and 1-covers $t$, so we assume that $z_t\in N_G(v)$. Since $t$
is adjacent from $v$  and $z_t$ is a problem for $t$, it follows that $z_t\notin N^+_G(v)$, so 
$z_t\in N^-_G(v)$.
Choose $y\in Y$ such that $y$ 2-covers 
$z_{t}$, and choose $u\in T$ such that $y\DD u\DD z_{t}$ is a directed path. Since $z_{t}$ is a problem for $t$, it follows
that $t$ is adjacent from $u$, and so $y$ 2-covers $t$. This proves (5).

\bigskip

Now let us bound the size of $K$. We have
$$|K|=1+ |Z\setminus D| + |S\setminus (Q\cup J\cup Z\cup D)|+ |X|+|Y|+|Q''\setminus D|.$$
We know that 
$$|F|+|S\setminus (Q\cup J\cup Z)|+ |Q''|+1\le |G|/2-|Z|,$$
and $|X|\le |F|$, and $|Y|\le |D\setminus (J\cup Q')|$.
Adding, we deduce that:
\begin{align*}
&|K|+ |F|+|S\setminus (Q\cup J\cup Z)|+ |Q''|+1 +|X|+|Y|\\
&\le 1+ |Z\setminus D| + |S\setminus (Q\cup J\cup Z\cup D)|+ |X|+|Y|+|Q''\setminus D|\\
&+ (|G|-|Z|-|B|)/2 + |F| + |D\setminus (J\cup Q')|.
\end{align*}
This simplifies to:
\begin{align*}
|K|+ |S\setminus (Q\cup J\cup Z)|+ |Q''|
&\le |Z\setminus D| + |S\setminus (Q\cup J\cup Z\cup D)|+ |Q''\setminus D|\\
&+ |G|/2-|Z| + |D\setminus (J\cup Q')|.
\end{align*}
Since $|Z|\le |B|$, and 
$$|S\setminus (Q\cup J\cup Z)| = |S\setminus (Q\cup J\cup Z\cup D)| + |D\setminus (Q\cup J\cup Z)|,$$
we deduce
$$|K|+ |D\setminus (Q\cup J\cup Z)|+ |Q''|
\le |Q''\setminus D|
 + |D\setminus (J\cup Q')|+ |Z\setminus D| + |G|/2-|Z|.$$
Since 
$$|D\setminus (J\cup Q')|-|D\setminus (Q\cup J\cup Z)|=  |(D\setminus J)\cap (Q''\cup (Z\setminus Q'))|\le |D\cap (Q''\cup Z)|,$$
this further simplifies to:
$$ |K|+ |Q''\cap D| \le |D\cap (Q''\cup Z)|- |Z\cap D|+ |G|/2,$$
and so
$ |K| \le |G|/2$.
This proves \ref{mainthm3}.~\bbox

\section{Large 2-kernels}\label{sec:large}

In this section, we turn to a second topic, Spiro's question \ref{largeconj}.
While it seems to be asking for something close to the opposite of \ref{smallconj}, Spiro observed that \ref{smallconj} implies 
\ref{largeconj}. Here is his argument: to prove \ref{largeconj} for a digraph $G$, choose a large number $n$. If $G$ has a source $v$, delete $v$ and all 
its out-neighbours and apply
induction; while if $G$ has no sources, for each vertex $v$ of $G$, add $n$ new vertices adjacent from $v$ and with no other neighbours.
Applying \ref{smallconj} with $n$ sufficiently large implies that $G$ satisfies \ref{largeconj}. 

If $G$ is a digraph and $X\subseteq V(G)$, let $N^+_G[X]$ denote the set of vertices that either belong to $X$
or are adjacent from a vertex in $X$.
The same construction (adding $nw(v)$
new out-leaves for each vertex) shows that \ref{smallconj} implies a slightly stronger statement ($\mathbb{Z}_+$ denotes the set of 
non-negative integers, and $f(X)$ denotes $\sum_{v\in X}f(v)$):
\begin{thm}\label{largeweightconj}
{\bf Conjecture:} In every digraph $G$, and for every map $f:V(G)\rightarrow \mathbb{Z}_+$ there is a 2-kernel $K$ such that $f(N^+_G[K])\ge f(V(G))/2$.
\end{thm}
In this section we show that \ref{largeweightconj} is true for split digraphs, and indeed for a somewhat more general class of graphs.
If $G$ is an oriented graph, let us say a {\em break} of $G$ is a partition $(S,T)$ of $V(G)$ such that $G[S]$ is {\em acyclic} 
(that is, has no directed cycles), and $G[T]$ is a tournament. We will show:
\begin{thm}\label{largethm}
In every oriented graph $G$ that admits a break, and for every map $f:V(G)\rightarrow \mathbb{Z}_+$, there is a 2-kernel $K$ such that 
$f(N^+_G[K])\ge f(V(G))/2$.
\end{thm}
The greater generality given by the function $f$ will be useful for the inductive proof, allowing us to delete vertices without changing $f(V(G))$.
We need a result of von Neumann and Morgenstern~\cite{neumann}:
\begin{thm}\label{neumann}
Every acyclic digraph has a unique 1-kernel.
\end{thm}
In order to prove \ref{largethm}, we prove a stronger statement (by the {\em non-neighbourhood} of a vertex $v$, we mean the digraph
induced on the set of vertices different from and nonadjacent with $v$):
\begin{thm}\label{betterlargethm}
Let $(S,T)$ be a break of an oriented graph $G$, and let $f:V(G)\rightarrow \mathbb{Z}_+$ be a map.
Then there is a 2-kernel $K$ such that $f(N^+_G[K])\ge f(V(G))/2$, where either $K\subseteq S$, or $K$ consists of some $v\in T$
together with the unique 1-kernel of its non-neighbourhood.
\end{thm}
\Proof
We assume the result holds for all oriented graphs that admit breaks $(S',T')$ with $2|S'|+|T'|< 2|S|+|T|$.
For each $X\subseteq S$, let $A(X)$ be the unique 1-kernel of $G[X]$ (which exists by \ref{neumann}); and for each $v\in T$, let $M(v)$ be 
its non-neighbourhood. Let us say a 2-kernel $K$ of $G$ is {\em special for $(G,S,T)$} if either
$K\subseteq S$, or $K=\{v\}\cup A(M(v))$ for some $v\in T$. 
\\
\\
(1) {\em We may assume that $\{v\}\cup A(M(v))$ is a 2-kernel for each $v\in T$.}
\\
\\
Suppose not. Certainly $\{v\}\cup A(M(v))$ is stable, so there is a vertex $w\ne v$ such that $\{v\}\cup A(M(v))$ does not 2-cover $w$.
We claim that $N^-_G(w)\subseteq N^-_G(v)$. For suppose that $s\in N^-_G(w)\setminus N^-_G(v)$. Since $s\notin \{v\}\cup N^+_G(v)$
(because $\{v\}\cup A(M(v))$ does not 2-cover $w$, it follows that $v,s$ are nonadjacent, and so $s\in M(v)\subseteq S$. But then $s$
is 1-covered by $A(M(v)$, and so $w$ is 2-covered by $\{v\}\cup A(M(v))$, a contradiction.
This proves that $N^-_G(w)\subseteq N^-_G(v)$.
Thus every 2-kernel of $G'=G\setminus v$ is also a 2-kernel of $G$. Define $f'(w) = f(w) + f(v)$,
and $f'(x) = f(x)$ for all $x\in V(G)\setminus \{v,w\}$. Applying the inductive hypothesis to $G'$ and $f'$,
we deduce there is a 2-kernel $K$ of $G'$ (and hence of $G$), special for $(G',S,T\setminus \{v\})$ (and hence special for $(G,S,T)$),
 such that $f'(N^+_{G'}[K])\ge f'(V(G))/2=f(G)/2$.
But $N^+_{G'}[K]\subseteq N^+_{G}[K]$, and if $w\in N^+_{G'}[K]$ then $v,w\in N^+_{G}[K]$, and so
$f'(N^+_{G'}[K])\le f(N^+_{G}[K])$. Hence $f(N^+_{G}[K])\ge f(G)/2$. 
This proves (1).

\bigskip
A {\em sink} of $G$ is a vertex that has no out-neighbours.
\\
\\
(2) {\em Let $s\in S$ be a sink of $G[S]$. We may assume that $s$ is a neighbour of every vertex in $T$.}
\\
\\
For each $t\in T$, if $s,t$ are nonadjacent, let us add the edge $ts$, forming an oriented graph $G'$. Suppose the
theorem holds for $G'$, with the same function $f$, and let $K'$ be a 2-kernel of $G'$, special for $(G',S,T)$, with 
$f(N^+_{G'}[K'])\ge f(V(G'))/2=f(V(G))/2$. For each $v\in T$, let $M'(v)$ be the non-neighbourhood of $v$ in $G'$.
There are four cases: 
\begin{itemize}
\item $K'=\{v\}\cup A(M'(v))$ for some $v\in T$ adjacent from $s$ in $G$;
\item $K'=\{v\}\cup A(M'(v))$ for some $v\in T$ adjacent to $s$ in $G$;
\item $K'=\{v\}\cup A(M'(v))$ for some $v\in T$ nonadjacent with $s$ in $G$;
\item $K'\subseteq S$.
\end{itemize}
In the first two cases, $M'(v)=M(v)$, and $\{v\}\cup A(M(v))$ is a 2-kernel of $G$ by (1); and $N^+_{G'}[K']=N^+_G[K']$, and so $K'$
satisfies the theorem. In the third case, $M'(v)=M(v)\setminus \{s\}$. If $A(M'(v))$ 1-covers $s$, then $A(M'(v))=A(M(v))$
and so $K'$ satisfies the theorem. If $A(M'(v))$ does not 1-cover $s$, then $A(M(v))=A(M'(v))\cup \{s\}$ (because $s$ is a sink of $G[S]$),
and so $K=\{v\}\cup A(M(v))$ satisfies the theorem. Finally, in the fourth case, $K'\subseteq S$. If $K'$ is a 2-kernel of $G$
then it satisfies the theorem, so we assume it is not; and since $K'$ is a 2-kernel of $G'$, it follows that
$K'$ does not 2-cover $s$. But then $K'\cup \{s\}$ satisfies the theorem. This proves (2).

\bigskip

If $S=\emptyset$, then $G$ is a tournament and the result holds, so we assume that $S\ne \emptyset$, and hence contains a sink of $G[S]$.
By (2), then $(S\setminus \{s\},T\cup \{s\})$ is also a break of $G$, and from the inductive hypothesis, 
there is a 2-kernel $K$ of $G$ such that $f(N^+_G[K])\ge f(V(G))/2$, and $K$ is special for $(G,S\setminus \{s\},T\cup \{s\})$.
But then $K$ is also special for $(G,S,T)$. This proves \ref{betterlargethm}.~\bbox

What happens to \ref{largeweightconj} if we assume that $V(G)$ can be partitioned into two sets $S,T$ where $T$ is a tournament
and $S$ is small? By \ref{betterlargethm},
the conjecture holds if $|S|\le 2$, and in hope of finding a counterexample, we worked on the case when $|S|=3$. But the conjecture is
also true in this case (by an {\em ad hoc} argument that does not seem capable of any generalization, and we omit the details).

There is a natural refinement of the conjectures \ref{smallconj} and \ref{largeconj}, equivalent to \ref{smallconj}
and implying \ref{largeconj}, that:
\begin{thm}\label{mixedconj}
{\bf Conjecture:} In every digraph $G$, and for every map $f:V(G)\rightarrow \mathbb{Z}_+$ there is a 2-kernel $K$ such that 
$|K|+f(V(G))/2\le |G|/2+ f(N^+_G(K))$.
\end{thm}
To deduce this from \ref{smallconj}, add $f(v)$ out-leaves to each vertex $v$. It implies \ref{smallconj} by taking $f(v)=0$ for all $v$,
and it implies \ref{largeconj} by scaling $f$ to be very large. Perhaps the proof of \ref{mainthm2} can be modified to show that
split graphs satisfy \ref{mixedconj}, but we have not seriously attempted this.

\section{$k$-kernels}\label{sec:kkernels}
Now we turn to the proof of our third result, \ref{arbthm}.
We begin with: 
\begin{thm}\label{dominator}
For all integers $k\ge 0$, if $G$ is an acyclic digraph with only one source, then there exists $X\subseteq V(G)$ with 
$|X|\le 1+(|G|-1)/(k+1)$ that $k$-covers $V(G)$.
Moreover, either $|G|=1$ or $|X|\le 1+(|G|-2)/(k+1)$ or $X$ is not stable.
\end{thm}
\Proof Let $r$ be the unique source. 
If $|G|\le k$, we may take $X=\{r\}$; then $|X|\le 1+(|G|-2)/(k+1)$ unless $|G|=1$, so the result holds. 
We assume then that $|G|>k$, and 
proceed by induction on $G$. For each $v\in V(G)$, let $A_v$ be the set of vertices that
are joined by a directed path (of any length) from $v$; and choose $v$ with $|A_v|$ minimal such that $|A_v|\ge k+1$.
(This is possible since $|A_r|\ge k+1$.) For each $w\in A_v$, there is a directed path $P$ from $v$ to $w$, and if $P$ has length 
more than $k$ then we may replace $v$ by its outneighbour in $P$, contradicting the minimality of $A_v$. Thus every vertex in $A_v$ is 
joined from $v$ by a path of length at most $k$. If $v=r$ then we may take $X=\{r\}$ and win as before, so we assume that $v\ne r$.
Let $G'$ be the digraph obtained by deleting $A_v$. Every vertex of $G'$ has an in-neighbour in $G'$ except $r$, so $G'$ has a 
unique source; and from the inductive hypothesis, there exists $X'\subseteq V(G')$ such that 
$|X'|\le 1+(|G'|-1)/(k+1)$ and $X'$ $k$-covers $V(G')$.
Moreover, either $|G'|=1$ or $|X'|\le 1+(|G'|-2)/(k+1)$ or $X'$ is not stable. Let $X=X'\cup \{v\}$. Thus $X$  $k$-covers $V(G)$.
Moreover, since $|A_v|\ge k+1$, it follows that
$|X|\le 1+(|G|-1)/(k+1)$, and if either $|X'|\le 1+(|G'|-2)/(k+1)$ or $X'$ is not stable, then correspondingly either
$|X|\le 1+(|G|-2)/(k+1)$ or $X$ is not stable. So we assume that $|G'|=1$, and so $V(G')=\{r\}$. 
Since $G$ has a unique source, it follows
that $v$ is adjacent from $r$, and so $X$ is not stable. This proves \ref{dominator}.~\bbox

We deduce:

\begin{thm}\label{acyclic2}
For every integer $k\ge 1$, if $G$ is an acyclic digraph with $|G|>1$ and with only one source, then $G$ has a $k$-kernel of 
size at most $1+(|G|-2)/k$.
\end{thm}
\Proof
By \ref{dominator} applied to $G$ with $k$ replaced by $k-1$, there exists $X\subseteq V(G)$ with 
$|X|\le 1+(|G|-1)/k$ that $(k-1)$-covers $V(G)$.
The digraph $G[X]$ is acyclic and hence has a 1-kernel $Y$, 
by \ref{neumann}.
Hence $Y$ is a $k$-kernel in $G$. Moreover, since $|G|\ge 2$, either
$|X|\le 1+(|G|-2)/k$ (when $|Y|\le |X|$ and the result is true), or $X$ is not stable (when 
$|Y|\le |X|-1\le (|G|-1)/k$ and again the result is true). This proves \ref{acyclic2}.~\bbox

As we said before, this result is tight (see figure \ref{fig:3kernels}).
Now let us deduce \ref{arbthm}, which we restate:
\begin{thm}\label{arbthm2}
For all integers $k\ge 2$, every digraph $G$ with $|G|>1$ and with a spanning arborescence has a $k$-kernel of size at most 
$1+(|G|-2)/(k-1)$.
\end{thm}
Since $G$ has a spanning arborescence, its vertex set can be numbered $\{v_1\LL v_n\}$ in such a way that for $2\le j\le n$
there exists $i\in \{1\LL j-1\}$ such that $v_iv_j$ is an edge. Let $A$ be the set of all edges $v_iv_j$ of $G$ with $i<j$, and let $B=E(G)\setminus A$.
Let $G_A$ be the subgraph with vertex set $V(G)$ and edge set $A$, and define $G_B$ similarly. Both $G_A,G_B$
are acyclic, and $G_A$ has a unique source. By \ref{acyclic2} applied to $G_A$ with $k$ replaced by $k-1$, 
$G_A$ has a $(k-1)$-kernel $X$ of size at most $1+(|G|-2)/(k-1)$. Now $X$ is stable in $G_A$,  and $G_B[X]$ is acyclic, 
and so has a 1-kernel $Y$, by \ref{neumann}. But then $Y$ is a $k$-kernel in $G$, and $|Y|\le |X|\le 1+(|G|-2)/(k-1)$. This proves \ref{arbthm2}.~\bbox

\end{document}